\documentclass[10pt,twoside]{article}
\usepackage{graphicx}
\usepackage{amsmath,amssymb}
\usepackage{Latex-document}

\title{\bf From Quaternions to Cosmology:\vskip -2mm\ Spaces of Constant
Curvature, \vskip -2mm ca. 1873--1925\vskip 6mm}

\author{Moritz Epple\thanks{Universit\"at Stuttgart, Germany. E-mail: epple@math.uni-bonn.de}\vspace*{-0.5cm}}

\date{\vspace{-8mm}}

\begin{document}

\maketitle

\thispagestyle{first} \setcounter{page}{935}

\begin{abstract}

\vskip3mm

After mathematicians and physicists had learned that the structure of physical space was not necessarily
Euclidean, it became conceivable that the global topological structure of space was non-trivial. In the context of
the late 19th century debates on physical space this speculation gave rise to the problem of classifying spaces of
constant curvature from a topological point of view. William Kingdon Clifford, Felix Klein and Wilhelm Killing,
the latter of whom devoted a substantial amount of work to the topic in the early 1890s, clearly perceived this
problem as relevant for both mathematics and natural philosophy (i.e., physics or cosmology). To some extent, a
cosmological interest may even be found among those authors who restated the space form problem in more modern
terms in the early 20th century, such as Heinz Hopf.

\vskip 4.5mm


\end{abstract}

\vskip 12mm

\section{Scientific contexts of topology}

\vskip-5mm \hspace{5mm}

The broader aim of the present paper is to contribute to a better understanding of the emergence of modern
topology. From its very beginnings, {\itshape analysis situs} or {\itshape Topologie}, as Johann Benedikt Listing
proposed to call the new field in the 1840s, was perceived as one of the most basic subfields of mathematics.
Conceptually independent of many other branches of mathematics, it deserved thorough research in its own right.
During the 20th century, this perception became even more pronounced with the gradual growth of structural
thinking in mathematics. As is well known, topology --- axiomatized in set-theoretical terms following the lead of
Felix Hausdorff --- became one of three ``mother structures'' in Bourbaki's architecture of mathematics, making
topology into a paradigm field of pure mathematics. Only in recent decades have the immediate connections of
topology with science, and physics in particular, been emphasized in many lines of research.

The purist view of topology also has dominated historical research on the emergence of topology for a long time.
However, historians have begun to move beyond a history of topology focusing exclusively or at least predominantly
on conceptual developments within pure mathematics. For instance, the importance of Poincar\'e's interest in
celestial mechanics for the development of his qualitative theory of differential equations and of a number of
crucial topological ideas (such as the notion of homoclinical points or his ``last geometric theorem'') has been
underlined in several historical studies (see, e.g., [1]). Another area where the interaction of topological
research and physics has been investigated thoroughly is the emergence of the theory of Lie groups [2]. Finally,
the relations between 19th century studies of vortex motion in ideal fluids by Helmholtz and Thomson, the latter's
influential theory of vortex atoms, and the early attempts at a classification of knots and links by Tait and his
followers have been studied in detail [3], [4].

This and similar research has made it clear that the gradual
formation of topology in the latter half of the 19th century and
the first decades of the 20th century involved more than just pure
mathematics. In addition to the growing need for topological
notions in fields such as (algebraic) function theory or
differential geometry, a need for topology was clearly felt in
several domains of physics (and maybe even in chemistry). In the
following, another major case will be discussed in which
mathematical and physical thinking jointly contributed to the
emergence of new topological ideas.\footnote{A more detailed
treatment of this episode, including full references, will appear
in [5].}

\section{The topological space problem}

\vskip-5mm \hspace{5mm}

In a well-known series of events ranging from the first
mathematical discussions of non-Euclidean geometries to heated
public debates in the late 19th century, mathematicians and
physicists learned that the most adequate mathematical description
of physical space was not necessarily Euclidean. This insight had
a wide range of consequences both for the body and for the image
of geometric, and indeed mathematical, knowledge (to use a
distinction proposed by Yehuda Elkana). One of these consequences
was to challenge not only the metric properties of Euclidean space
(as a model of physical space) but to question its other
properties as well. If there was no {\itshape a priori} reason for
accepting the axiom of parallels, why should there be {\itshape a
priori} reasons for accepting, e.g., the topological features of
Euclidean space? What were the topological types of the best
mathematical descriptions of physical space?

To phrase such a question in modern terms and within our
understanding of the relations between geometry and topology
sounds anachronistic. Nevertheless, a corresponding problem
{\itshape was} raised in the terms available to 19th century
scientists. Before the establishment of a coherent framework of
topological notions, such terms were, in particular, the
dimension, the  ``Zusammenhang'' or connectivity, and the
continuity of space. While it is well known that the issue of the
dimension of space was in the focus of several 19th century
debates, it has less often  been emphasized that the properties of
connectivity and continuity of space came into question as well.
Here I will concentrate on the problem of the connectivity of
space.

The notion of ``Zusammenhang'' was originally introduced by
Bernhard Riemann as a tool for distinguishing different types of
(Riemann) surfaces in the context of function theory. This notion
does not figure prominently in his famous talk {\itshape \"Uber
die Hypothesen, welche der Geometrie zu Grunde liegen}. There,
Riemann introduced the crucial distinction between the
``Ausdehnungsverh\"altnisse'' and ``Ma{\ss}verh\"altnisse''
(roughly:  topological properties vs. metric properties) of a
manifold, but the global topological aspects of manifolds received
no special emphasis. This holds in particular for the final
sections of his talk which were devoted the geometry of physical
space. While acknowledging that there exists a ``discrete
manifold'' of possible ``Ausdehnungsverh\"altnisse'' of space,
Riemann expressed scepticism about pursuing the global properties
of space beyond the issue of dimension: ``Questions about the
immeasurably large are idle questions for the explanation of
nature.''

When Riemann's talk reached the scientific public in 1868, another
contribution that shaped the later debates on the space problem
was on its way. The physicist Hermann v. Helmholtz argued that for
epistemological reasons, a crucial assumption in any mathematical
description of physical space should be the ``free mobility of
rigid bodies'' of arbitrary size. According to Helmholtz, the
existence of freely movable rigid bodies was a precondition for
measuring lengths. In mathematical terms, it implied that the
classical non-Euclidean geometries were the only possible models
of space. Although Helmholtz's argument was soon criticised for
technical reasons, his main assumption (not easily stated in
precise mathematical terms) was accepted during the 19th century
even by many proponents of liberal approaches to the geometry of
physical space. With one exception and one crucial modification,
this holds for all authors that will be treated below.

The exception is William Kingdon Clifford, the most imaginative
follower of Riemann's geometric speculations in Britain. His
remarks on a ``space theory of matter'', according to which all
material phenomena might be explained by a time-dependent,
wave-like variation of space curvature, are well known. In
addition, several of Clifford's writings show a marked interest in
different {\itshape global} possibilities for manifolds or spaces.
In [6], Clifford hinted at a large variety of ``algebraic
spaces'', higher dimensional analogues of Riemann's surfaces. In
the same paper, he presented his example of a closed surface
embedded in elliptical 3-space, the inherited geometry of which is
locally Euclidean. As this example came to play an essential role
in the following, let me recall the main line of Clifford's
construction.

Identifying points in elliptic space with with one-dimensional
subspaces of the quaternions, any given quaternion different from
zero induces two isometries of elliptic space by left and right
multiplication. Such isometries Clifford called left and right
``twists'', respectively. (Felix Klein would later term them
``Schiebungen'', translations.) Every twist possessed a
space-filling family of invariant lines, i.e. it moved points
along these invariant lines by a constant distance. Any two
members of one and the same such family were called ``parallels''
by Clifford. Next, given any two intersecting lines $l$ and $l'$,
Clifford considered the ruled surface generated by all those
Clifford (left) parallels to $l$ which met $l'$. Equivalently,
this surface could be described as being generated by all (right)
Clifford parallels to $l'$ meeting $l$. Moreover, there were two
commuting one-parameter families of left and right twists inducing
isometries of the surface, which had $l$ and $l'$ as invariant
lines, respectively. Consequently, the surface had constant
curvature zero. In topological terms, the surface was a torus as
may be seen from Clifford's description of it as ``a finite
parallelogram whose opposite sides [given by the lines $l$ and
$l'$] are regarded as identical'' [6, p.  193]. Closer inspection
shows that the surface is indeed orientable. It is important to
keep in mind that Clifford's example was not constructed by
endowing the 2-torus with a geometrical structure, but rather as a
particular surface embedded in elliptic 3-space arising from the
consideration of a particular set of isometries, Clifford's twists
or translations.

Several remarks in Clifford's philosophical articles indicate that
he was aware of the implications this example had for the problem
of giving an adequate mathematical description of physical space:
The same local geometry might be tied to spaces that are globally
different. Even for spaces of constant curvature one could make
different ``assumptions [...] about the {\it Zusammenhang} of
space'', as he wrote in 1873 [7, p. 387]. Clifford also saw that
these differences were of a topological nature. A remark of 1875
may even be read as advocating a more radical kind of
`topologism': ``There are many lines of mathematical thought which
indicate that distance or quantity may come to be expressed in
terms of {\it position} in the wide sense of the {\it analysis
situs}. And the theory of space-curvature hints at a possibility
of describing matter and motion in terms of extension only.'' [7,
p. 289.]

In the 1870s, Cliffords critical remarks about the possibilities of globally different spaces with the same local
geometry seem not to have generated resonances within the scientific communities either of physicists or of
mathematicians. This changed during the 1880s for reasons that originally had nothing to do with Clifford's ideas.
In 1877, the American astronomer Simon Newcomb published a paper on a geometry of space with constant positive
curvature (in Kleinian terms: elliptic geometry). In a reaction to this paper, Wilhelm Killing, a student of
Weierstrass and mathematics teacher, argued that Newcomb had overlooked the fact that there were actually
{\itshape two} possible geometries with constant positive curvature that should be discussed: elliptic and
spherical space. This prompted Felix Klein, whose earlier contributions on non-Euclidean geometry also had focused
on elliptic rather than spherical geometry, to enter into a correspondence with Killing.\footnote{Killing's
letters to Klein may be found in the Nieders\"achsische Staats- und Universit\"atsbibliothek (NSUB) G\"ottingen,
Handschriftenabteilung, Cod. MS Klein 10.} While Klein pointed out that Killing's remark was fairly obvious from
the perspective that Klein had developed, Killing repeatedly emphasized the importance of a theorem (along
Helmholtz's line of argument) specifying the full range of geometric spaces compatible with the idea of the free
mobility of rigid bodies. According to Killing, there were exactly four such spaces: 3-dimensional Euclidean,
hyperbolic, elliptic, and spherical space. Killing was clearly interested in what might be called the foundations
of physical geometry as opposed to the framework of projective geometry that guided Klein. One of Klein's
reactions now was to refer to Clifford's flat surface in elliptic space. In his eyes, this example showed that
there were many more manifolds satisfying the assumptions that Killing wanted to hold. Killing protested:
Clifford's surface did not admit free mobility in the full sense (it did not allow global rotations) and thus was
not a ``space form satisfying our experience'' (Killing to Klein, cf. note 2, 5 October 1880).

It took Killing and Klein several years to sort out their differences. In the end it became clear (not least
because of Sophus Lie's additional work on Helmholtz's approach) that the conditions of constant curvature and
free mobility in the Helm\-holtzian sense had to be distinguished. The former was a local, the latter both a local
{\itshape and} a global property of space. However, Klein pointed out that there was no clear empirical sense
which could be given to this latter property --- contrary to both Helmholtz's and Killing's intentions. What
{\itshape might} make sense as an empirical requirement was the free mobility of bodies of finite size, indeed of
globally bounded finite size. Of course this restricted condition of free mobility still implied a constant
curvature of space. Hence Klein felt justified in posing the following problem, first in a lecture course on
non-Euclidean geometry in 1889/1890, then in print: ``to enumerate all species of connectivity which may at all
occur in closed manifolds of some constant measure of curvature'' [8, p. 554]. Obviously, Klein was interested in
the global topological differences of such manifolds, not in a finer classification up to isometry. In his paper,
he gave a (not quite complete) discussion of the two-dimensional case, emphasizing again Clifford's work. Then he
pointed out the general connection between regular tessellations of the standard non-Euclidean spaces of dimension
3 and manifolds of constant curvature. From his own and from Poincar\'e's work on automorphic functions he knew
that this connection lead to quite involved problems. The corresponding section of his paper included an
invitation ``that the question would be taken up elsewhere''. He underlined that the problem was ``fundamental for
the doctrine of space, inasmuch as we want to start the latter from the condition of free mobility of rigid
bodies'' [8, p. 564].

Note that Klein here refered to the {\itshape restricted}
condition of free mobility. By now, Killing accepted Klein's
argument that only this latter version of the condition had
empirical content, and in the following years he took up the task
that Klein had set. One may group his work on what he now called
the problem of ``Clifford-Klein space forms'' (in the following:
CK space forms) under three headings: a reformulation of the
problem in group-theoretical terms, the construction of new
classes of examples, and a discussion of the scientific relevance
of spaces of constant curvature. I will return to the two more
mathematical aspects in the next section. Here I want to comment
on the third.

In both of his relevant publications, Killing included long
sections defending a study of CK space forms in the context of the
foundations of physical science [9], [10, part 4]. Repeating
Klein's argument, Killing advocated an understanding of free
mobility in the restricted sense and emphasized that nothing in
experience excluded the possibility of space being different from
the standard non-Euclidean spaces. In fact, he considered only one
possible criticism as requiring a more careful discussion: As yet,
neither mechanics nor any other physical theory existed for CK
space forms. But this was just the usual course of science. For
the standard non-Euclidean spaces as well, mechanics was just in
the process of being developed (Killing himself had made important
contributions). In consequence, the primary task was to develop
physical theories for CK space forms as well. Only then it would
be possible to judge their scientific merits. As a particular
phenomenon that mechanics in multiply connected CK space forms
might bring up, Killing mentioned anisotropies of the
gravitational force between two bodies [10, p. 347]. One may well
read this as a hint at a possible local empirical phenomenon that
might help in finding out global `connectivity properties' of
physical space.

Killing was not alone in the 1890's in discussing the physical
relevance of spaces of constant curvature. In 1899, Klein came
into contact with the young and aspiring astronomer Karl
Schwarzschild when the latter gave a talk at a large meeting of
astronomers discussing ``the admissible measure of curvature of
space''. In this talk, Schwarzschild gave bounds on the radii of
curvature of either an elliptic or a hyperbolical universe
consistent with astronomical observations of star parallaxes.
After the talk and in ensuing correspondence\footnote{See
Schwarzschild's letters in NSUB G\"ottingen, Cod. MS Teubner 44
and Cod. MS Klein 11.}, Klein made Schwarzschild aware of the fact
that in such a discussion, CK space forms should also be taken
into account. Schwarzschild agreed. In the printed version of his
talk, he added an appendix in which he briefly discussed whether
or not space might actually be a non-standard space of constant
curvature. In very intuitive terms, he explained to his readers
(the paper was published in an astronomical journal) how one could
conceive of astronomical observations suggesting such kinds of
spaces: by observing ``identical, apparent repetitions of the same
world-whole, be it in a Euclidean, elliptic, or hyperbolic
space''. However, the time was not yet ripe for a full discussion
of this possibility: ``We may treat the other Clifford-Klein space
forms very briefly, the more so since they have not yet been
investigated completely even from a mathematical point of view.
[...] experience only imposes, in all cases, the condition that
their volume has to be larger than that of the visible star
system.'' [11, appendix.]

Killing's mathematical work on CK space forms was reformulated in
modern mathematical terms and substantially extended by Heinz Hopf
who devoted one of two parts of his dissertation to the problem in
1925. Again I defer a discussion of the mathematical parts of
Hopf's work to the next section. However, it must be pointed out
that Hopf also shared a cosmological interest in CK space forms
with Killing and Klein. When, in 1928, Klein's lectures on
non-Euclidean geometry were edited posthumously in a completely
rewritten form by Walter Rosemann, Hopf took over the task of
writing a new section on spaces of constant curvature (at the time
called ``homogenous spaces'' by him). He closed this section with
discussing ``the application of geometry to the external world''.
Here, only ``the possibility of homogenous space forms [had] to be
taken into consideration'', as no empirical data were known that
would force one to consider spaces of variable curvature. Of
particular value was the ``possibility of ascribing to the
universe a finite volume, independently of its geometrical
structure [...] since the idea of an infinite extent [...] causes
various difficulties, for instance in the problem of the
distribution of mass.'' [12, p. 270.] One should note that this
was written after the advent of Einstein's theory of general
relativity, and after the development of relativistic cosmology
had seriously begun with contributions by Einstein, Schwarzschild,
de Sitter, Weyl and others. In this context, constant curvature
was no longer a pre-condition of measurement, but rather a
consequence of the assumption of a homogenous average distribution
of mass throughout the universe.

\section{Killing's and Hopf's mathematical contributions}

\vskip-5mm \hspace{5mm}

Killing's main mathematical contribution to the problem of CK
space forms was its reduction to group theoretical terms. Killing
tried to show that in all dimensions $n$, Klein's problem (see
above) was equivalent to finding all finitely or at most countably
generated subgroups $G$ of $SL(\mathbb R, n+1)$ which for some
real parameter $1/k^2$ leave invariant the bilinear form
\[ a(x,y) = k^2x_0y_0 + x_1y_1 + ... + x_ny_n\:, \qquad x,y\in
\mathbb R^{n+1}\:,\] and satisfy a discontinuity condition that
will become clear as we go along [10, p. 322]. Even if it is
difficult to follow all details of Killing's argument, its main
line is clear. (In the following, modern abbreviations are used to
condense Killing's verbal style.) To begin with, if $M$ was a
manifold of dimension $n$ with constant curvature $1/k^2$, Killing
required the existence of some $r>0$ such that for all points
$P\in M$ there existed a ball $B_r (P) \subset M$ of radius $r$
isometric with a similar ball in Euclidean, hyperbolic or
spherical space of the same curvature. This was Killing's way of
stating the restricted condition of free mobility. It implied both
a kind of completeness of the manifold and the discontinuity
condition just mentioned. Using a technique he had learned in a
seminar of Weierstrass, Killing translated this into local
``coordinates'', by which he understood isometric mappings
\[ B_r(P) \longrightarrow X_k := \{\: x \in \mathbb
R^{n+1} \:|\: a(x,x) = k^2 \:\}
\]
mapping $P$ to $\bar P := (1, 0, ..., 0)$. (For negative
curvature, the condition $x_0>0$ was added in the definition of
$X_k$; in the flat case, Killing just considered the hyperplane in
$\mathbb R^{n+1}$ defined by $x_0=1$.) Endowing $X_k$ with the
metric $d$ given by
\[ k^2 \cos \frac{d(x,y)}{k} = a(x,y)\:, \qquad x,y\in
X_k\:, \] made $X_k$ into a model of the standard Euclidean and
non-Euclidean spaces that had been used by Killing in most of his
earlier work on these geometries.

Choosing a particular point $P$ as the origin, these
``Weierstrassian coordinates'' defined a 1-1 correspondence of the
bundles of geodesics through $P$ and $\bar P$. Using this
intuition, Killing extended a local coordinate system around $P$
to a kind of global coordinate system, i.e., a `mapping' $M
\rightarrow X_k$, associating a point $Q\in M$ on some geodesic
through $P$ with a point $\bar Q \in X_k$ on the corresponding
geodesic such that the distances between $P$ and $Q$ and between
$\bar P$ and $\bar Q$ were equal. Partly without further argument
and partly based on intuitive explanations, Killing assumed that
this `mapping' was in general multi-valued (since in $M$ there
might exist closed geodesics), surjective, and locally isometric.
If $\bar Q_1, \bar Q_2 \in X_k$ were two ``coordinates'' of the
same $Q\in M$, then by construction there existed a (local)
isometry $K_r (\bar Q_1) \rightarrow K_r (\bar Q_2)$. Killing
assumed that this mapping could be uniquely extended to a global
isometry $\psi: X_k \rightarrow X_k$. He knew that isometries of
$X_k$ were induced by elements of the group we denote by
$SL(\mathbb R, n+1)$, leaving the bilinear form $a$ invariant.
Again on intuitive grounds Killing argued that any such $\psi$ was
in fact what we would call a covering transformation of $M$. The
collection of all $\psi$ arising in this way formed a discrete
subgroup $\Gamma$ of the isometry group of $X_k$ which had the
property that every $\psi \in \Gamma$ moved points by a distance
of at least $r$. $M$ itself was then equivalent to what later was
called the quotient space $X_k/\Gamma$.

The gaps and intuitive turns in Killing's argument give a striking
illustration of the growing need for precise topological arguments
in some areas of mathematics at this time. Notions relating to
covering spaces or the fundamental group (in his intuitive
explanations, Killing repeatedly relied on the consideration of
``motions of bodies'' along closed geodesics in $M$) would have
helped Killing significantly in securing the vaguer parts of his
considerations.

Killing was quite clear that the new problem in group theory was
difficult. Accordingly, he was satisfied with describing a few
simple Euclidean and spherical space forms. In the flat case, his
main example was the analogue of Clifford's surface, the manifold
given by $\mathbb R^3/\mathbb Z^3$. In the case of positive
curvature, Killing noticed that in even dimensions, only $S^n$ und
$\mathbb RP^n$ with their canonical metrics could occur. In
dimension 3 he mentioned other  possibilities, e.g. $\mathbb
RP^3/\Gamma$, where  $\Gamma$ is a cyclic group of Clifford's
translations.

It was Heinz Hopf who reworked Killing's arguments in a modern
framework. In his dissertation of 1925, he presented a completely
revised  treatment of the problem of CK space forms that made it
superfluous to look into the older literature any more [13]. His
version of the problem was to classify all geodetically complete
Riemannian manifolds of constant curvature in either of two
possible senses: One could try to classify the resulting
``geometries'' (i.e., look for a classification up to isometry) or
one might wish to classify just the manifolds carrying these
geometries (i.e., look for a classification up to diffeomorphism).
In Killing's work, this distinction had never been clearly made.

After a preliminary clarification of the relation between
Killing's earlier completeness condition and the weaker condition
of geodetic completeness, Hopf gave a new proof of Killing's basic
result. In the new setting, this theorem took the form that every
geodetically complete Riemannian manifold of constant curvature
(again called CK space form by Hopf) was a quotient of Euclidean,
hyperbolic or spherical space by a discontinuous group $\Gamma$ of
isometries without fixed points and such that no orbit of $\Gamma$
had a limit point. The geometric content of Hopf's proof was very
similar to Killing's argument -- the difference being that Hopf
had conceptual tools at his disposal that Killing had missed. Hopf
showed that if $M$ was a CK space form in his sense, then every
point $P\in M$ still had a neighbourhood that could be mapped
isometrically onto a neighbourhood of some point $\bar P$ in one
of the standard spaces, say $X$. Using again the resulting 1-1
correspondence of the bundles of geodesics through $P$ and $\bar
P$, Hopf defined a mapping $X\rightarrow M$ (!) of which he showed
that it was an isometric covering. As $X$ was simply connected, it
was the universal covering space of $M$. Moreover, the fundamental
group $\pi_1(M)$ acted freely and discontinuously in the sense
explained above by isometries on $X$. Therefore, $M$ was isometric
to a quotient manifold of the required form.

Instead of multi-valued ``coordinates'', Hopf could speak of
coverings, inverting the direction of the crucial mapping. Relying
on the notions of universal covering spaces, covering
transformations, local isometries etc., Hopf was able to formulate
several steps in his proof as simple arguments by contradiction.
For this step into mathematical modernity, Hopf had an important
model: Hermann Weyl. The conceptual framework used by Hopf was
mainly an adaptation of that outlined in the topological sections
of Weyl's monograph {\itshape Die Idee der Riemannschen Fl\"ache}
of 1913. On the other hand, Weyl himself had pointed out that this
framework could be used in discussing manifolds of constant
curvature. Using the language of coverings, Weyl showed in an
appendix to [14] that every ``closed Euclidean space'' (i.e.,
every closed manifold of constant curvature zero) was isometric to
a ``crystal'', i.e., a quotient of Euclidean space by a suitable
discrete group of isometries. Weyl's proof only worked in the flat
case, and apparently he did not pursue analogous questions for
curved manifolds.

Hopf not only used modern topological tools for reformulating the
space form problem. He also showed how topology could profit from
the geometric ideas involved in this problem. Looked at in this
perspective, one point was fairly obvious: CK space forms,
constructed as quotient spaces, furnished new examples of
manifolds with known fundamental groups. At the time, this was
particularly interesting for finite groups; very few examples of
manifolds with finite fundamental groups had been known beyond
those with cyclic groups. Based on Klein's analysis of the
isometry group of elliptic 3-space (which in turn relied on
Clifford's ideas), Hopf discussed a series of new 3-dimensional,
spherical space forms with finite fundamental groups as well as
infinitely many spaces with infinite groups.\footnote{Spherical,
3-dimensional space forms with polyhedral fundamental groups had
already been added to Killing's examples by the American geometer
Woods in 1905.} Moreover, Hopf's research on spaces of constant
positive curvature proved to be of decisive importance at a later
point in his career: Clifford's parallels provided him with the
crucial example of a mapping from $S^3$ to $S^2$ that was not
homotopic to a constant mapping (lifting a fibration of elliptic
3-space by Clifford parallels to the 2-sheeted covering of
elliptic 3-space produced what became known as the Hopf fibration
of $S^3$). Even Hopf's basic invariant for classifying such maps,
the linking number of fibres, was derived from an intuitive
understanding that the linking behaviour of Clifford's parallels
was an obstacle for deforming the Hopf fibration into a constant
[16].

\section{Conclusions}

\vskip-5mm \hspace{5mm}

In these ways, an important strand in the formation of modern
manifold topology profited from geometric ideas that had their
origin in a 19th century context in which mathematical and
cosmological thinking were closely related. If I may use a
topological metaphor: Once more it turns out that the fibres of
historical developments in different scientific fields are
intertwined rather than grouped in a locally trivial bundle, the
base of which would be some eternal architecture of concepts or
structures that only need straightforward elaboration.

From cases like the one I have sketched one might learn that in
periods of innovative research, boundaries between different
mathematical fields and even between mathematical and physical
thinking may tend to blur. In other periods, this may be
different. In the present case, it seems (at least at first sight)
that the later {\itshape solutions} of the Euclidean and the
spherical space form problems were found in episodes of
autonomous, purely mathematical research. Moreover, the history of
the space form problem between Clifford and Hopf reveals a
complicated relation between tradition and modernization within
mathematics. The analysis of Killing's and Hopf's ways of
approaching the space form problem shows that despite the crucial
differences due to the non-availability or availability of precise
topological notions, Killing's more traditional geometric ideas
were taken up in the modern formulations. Another such core of
geometric ideas that was handed down to modern topology were
Clifford's geometric ideas: in fact {\itshape all} later authors
discussed here made some use or other of these ideas.

In a period in which the global topological properties of the
universe receive new interest among cosmologists, it seems fitting
to recall that a century ago, such an interest also motivated some
of the topological problems and ideas mathematicians have since
become acquainted with. At least in this indirect way, questions
about the immeasurably large have {\itshape not} been idle
questions for the explanation of nature.

\label{lastpage}

\end{document}